\documentclass{article}
\usepackage{arxiv}

\usepackage{graphicx}
\usepackage{amsfonts}
\usepackage{amsmath}
\usepackage{amsthm}
\usepackage{amssymb}           
\usepackage{color}             
\usepackage{mathtools}         
\usepackage{booktabs}          
\usepackage{multirow}          
\usepackage{algorithm}
\usepackage{algpseudocode}
\usepackage{hyperref}          
\usepackage{setspace}          

\newcommand{\Mini}{\mathop{\rm Minimize}}

\newcommand{\grad}{\mathop{\rm grad}}
\newcommand{\Hess}{\mathop{\rm Hess}}

\newcommand{\D}{\mathrm{D}}
\DeclareMathOperator{\id}{id}
\DeclareMathOperator{\St}{St}
\DeclareMathOperator{\Sym}{Sym}
\DeclareMathOperator{\sym}{sym}
\DeclareMathOperator{\tr}{tr}
\DeclareMathOperator{\vect}{vec}
\DeclareMathOperator{\qf}{qf}

\newtheorem{assumption}{Assumption}
\newtheorem{definition}{Definition}
\newtheorem{remark}{Remark}
\newtheorem{proposition}{Proposition}

\newtheorem{lemma}{Lemma}
\newtheorem{theorem}{Theorem}

\title{Modified Armijo line search in optimization on Riemannian submanifolds with reduced computational cost
}

\author{Hiroyuki Sato \\
Department of Mathematical Sciences\\
Ritsumeikan University\\
1-1-1 Noji-higashi, Kusatsu-shi, Shiga 525-8577, Japan \\
\texttt{hsato@fc.ritsumei.ac.jp} \\
	\And
{Yuya Yamakawa} \\
Graduate School of Management\\
Tokyo Metropolitan University\\
1-1 Minami-Osawa, Hachioji-shi, Tokyo 192-0397, Japan\\
\texttt{yuya@tmu.ac.jp} \\
\AND
Kensuke Aihara \\
Department of Computer Science\\
Tokyo City University \\
1-28-1 Tamazutsumi, Setagaya-ku, Tokyo 158-8557, Japan\\
\texttt{aiharak@tcu.ac.jp} \\
}

\date{\today}

\renewcommand{\headeright}{}
\renewcommand{\undertitle}{}
\renewcommand{\shorttitle}{Modified Armijo line search in optimization on Riemannian submanifolds with reduced computational cost}

\begin{document}

\maketitle

\begin{abstract}
For optimization problems on Riemannian manifolds, many types of globally convergent algorithms have been proposed, and they are often equipped with the Riemannian version of the Armijo line search for global convergence.
Such existing methods need to compute the value of a retraction mapping regarding the search direction several times at each iteration; this may result in high computational costs, particularly if computing the value of the retraction is expensive.
To address this issue, this study focuses on Riemannian submanifolds of the Euclidean spaces and proposes a novel Riemannian line search that achieves lower computational cost by incorporating a new strategy that computes the retraction only when inevitable.
A class of Riemannian optimization algorithms, including the steepest descent and Newton methods, with the new line search strategy is proposed and proved to be globally convergent.
Furthermore, numerical experiments on solving optimization problems on several types of Riemannian submanifolds illustrate that the proposed methods are superior to the standard Riemannian Armijo line search-based methods.
\end{abstract}

\keywords{Riemannian optimization \and steepest descent method \and Newton method \and Armijo line search \and low computational cost \and global convergence}

\section{Introduction} \label{sec:intro}
This study focuses on solving unconstrained optimization problems on Riemannian manifolds embedded in Euclidean spaces.
Since the 2000s, many researchers have studied optimization on Riemannian manifolds and have extended well-established methods for unconstrained optimization in Euclidean spaces, including both deterministic and stochastic ones, to those on Riemannian manifolds~\cite{AbsMahSep2008,boumal2023introduction,sato2021riemannian,sato2019riemannian}.
Such extended optimization on Riemannian manifolds is known as {\it Riemannian optimization}, and these techniques have been implemented in several software programs~\cite{bergmann2022manopt,boumal2014manopt,huang2018roptlib,martin2016manifoldoptim,townsend2016pymanopt}.

In this paper, we focus on line search-based iterative optimization methods on Riemannian manifolds.
In such methods, a sequence of iterates $\{x_k\}$ is generated on a Riemannian manifold ${\cal M}$ and two elementary procedures, computing a search direction $p_{k}$ belonging to a tangent space at $x_{k}$ and finding a step length $\alpha_k > 0$, are performed.
Many approaches, such as the steepest descent~\cite{AbsMahSep2008,boumal2023introduction,sato2021riemannian}, conjugate gradient~\cite{sato2021riemannian,sato2022riemannian,zhu2020riemannian}, and Newton's methods~\cite{AbsMahSep2008,boumal2023introduction}, have been proposed to compute specific search directions in optimization.
After determining a search direction, we find an appropriate step length, which is the primary focus of this paper.

Given a point $x_k$ on ${\cal M}$ and search direction $p_{k}$, a retraction $R$ (see Definition~\ref{def:retraction}) is used to define a curve $\gamma_k(\alpha) \coloneqq R_{x_k}(\alpha p_k)$ ($\alpha \in \mathbb{R}$) on ${\cal M}$, which satisfies $\gamma_k(0) = x_k$ and $\dot \gamma_k(0) = p_k$.
On a curve $\gamma_k$, when an appropriate point $\gamma_{k}(\alpha_k) = R_{x_k}(\alpha_kp_k)$ for some step length $\alpha_k > 0$ is found, it is adopted as the next point $x_{k+1} \in {\cal M}$.
Therefore, an appropriate step length $\alpha_k > 0$ should be efficiently determined when searching on the curve $\gamma_k$.
To determine whether a step length is appropriate, certain conditions are imposed on step length $\alpha_k$; for example, $\alpha_k$ is such that the objective function value $f(R_{x_k}(\alpha_k p_k))$ is sufficiently smaller than $f(x_k)$.
Among such conditions, the Armijo condition~\cite{AbsMahSep2008,boumal2023introduction,sato2021riemannian} is simple but crucial.
To determine $\alpha_k$ satisfying the Armijo condition, we typically use the backtracking approach---with an initial guess $\bar{\alpha} > 0$ for $\alpha_k$, we repeatedly contract the current guess by a contract factor $\beta \in (0, 1)$.
Thus, we test step lengths $\bar{\alpha}, \beta \bar{\alpha}, \beta^2 \bar{\alpha}, \dots$ until we ascertain that the instance $\beta^{l_k} \bar{\alpha}$ satisfies the Armijo condition for some nonnegative integer $l_k$.

To determine if a trial step length $\beta^l \bar{\alpha}$ satisfies the Armijo condition, the objective function value $f(R_{x_k}(\beta^l \bar{\alpha}))$ should be evaluated at the trial point $R_{x_k}(\beta^l \bar{\alpha})$.
Therefore, before $\beta^{l_k} \bar{\alpha}$ is adopted as the step length $\alpha_k$ to be used for computing $x_{k+1} = R_{x_k}(\alpha_k p_k) = R_{x_k}(\beta^{l_k}\bar{\alpha}p_k)$, trial points $R_{x_k}(\bar{\alpha} p_k)$, $R_{x_k}(\beta\bar{\alpha} p_k)$, \dots, $R_{x_k}(\beta^{l_k-1}\bar{\alpha} p_k)$ need to be computed and then discarded.
In the Euclidean case, that is, when ${\cal M} = \mathbb{R}^n$, a natural retraction is defined for $x, p \in \mathbb{R}^n$ as $R_{x}(p) \coloneqq x + p$, and computing trial points $x+\bar{\alpha} p$, $x+\beta \bar{\alpha} p$, \dots, $x+\beta^{l_k-1} \bar{\alpha} p$ is not expensive.
By contrast, computing the trial points $R_{x_k}(\beta^l \bar{\alpha} p_k)$ for $l = 0, 1, \dots, l_k-1$ with a retraction $R$ on a manifold is generally expensive.
Therefore, in this paper, we propose a novel approach for reasonably reducing this cost by computing a retraction only when inevitable.

Note that, in the aforementioned standard Riemannian optimization framework, the retracted points $R_{x_k}(\beta^l \bar{\alpha} p_k)$ for $l = 0, 1, \dots, l_k - 1$ are finally discarded and not used for the subsequent process.
Therefore, we posit that such points can be roughly computed as long as a step length satisfying the Armijo condition can be found.
Here, assuming that the manifold ${\cal M}$ in consideration is an embedded manifold of the Euclidean space $\mathbb{R}^n$, both a point $x_k$ on ${\cal M}$ and its tangent vector $p_k$ can be regarded as vectors in $\mathbb{R}^n$; therefore, $x_k + \beta^l \bar{\alpha} p_k$ can be considered as a vector in the ambient space $\mathbb{R}^n$, which can be used to approximate $R_{x_k}(\beta^l \bar{\alpha} p_k)$.
Furthermore, computing $x_k + \beta^l \bar{\alpha} p_k$ is generally considerably cheaper than computing $R_{x_k}(\beta^l \bar{\alpha} p_k)$.
The proposed method is based on this concept.

In summary, instead of $R_{x_k}(\beta^l \bar{\alpha} p_k)$, the proposed method uses $x_k + \beta^l \bar{\alpha} p_k$  and the corresponding approximate version of the standard Armijo condition, which avoids the computation of a retraction.
However, a step length satisfying the standard Riemannian Armijo condition is eventually required.
Therefore, after determining a step length $\tilde{\alpha}_k$ that satisfies the approximate Armijo condition, we compute $R_{x_k}(\tilde{\alpha}_k p_k)$ and check whether $\tilde{\alpha}_k$ satisfies the standard Riemannian Armijo condition at $x_k$.
If it is satisfied, we set $\alpha_k \coloneqq \tilde{\alpha}_k$ and proceed with the subsequent process using this step length; otherwise, we again seek a better trial step length using the approximate Armijo condition.

The remainder of this paper is organized as follows.
In Section~\ref{sec:preliminaries}, we introduce the notation and terminology used in this paper.
Some Riemannian optimization concepts are also reviewed.
In Section~\ref{sec:Armijo}, the framework of the backtracking-based Armijo line search algorithm on a manifold and its drawbacks are reviewed.
An improved strategy is then proposed.
In particular, we present a class of optimization methods, including the steepest descent, Newton, and quasi-Newton methods based on our efficient Armijo line search.
The global convergence of the proposed method is analyzed in Section~\ref{sec:Convergence}.
In Section~\ref{sec:experiements}, we present numerical experiments, in which the proposed methods are compared with the existing methods by applying them to Riemannian optimization problems on the sphere, the Stiefel manifold, and the manifold of symmetric positive definite matrices.
Section~\ref{sec:Conclusion} concludes the paper.

\section{Preliminaries}
\label{sec:preliminaries}
In this section, we first introduce the notation and terminology related to Riemannian optimization. Second, we define an unconstrained optimization problem on Riemannian manifolds and provide a prototype algorithm to solve the optimization problem.

\subsection{Notation and terminology}
The set of natural numbers (positive integers) is represented by $\mathbb{N}$.
For $n \in \mathbb{N}$, we regard ${\mathbb{R}^n}$ as an $n$-dimensional Euclidean space equipped with the standard inner product $\langle \cdot, \cdot \rangle \colon {\mathbb{R}^n} \times {\mathbb{R}^n} \to \mathbb{R}$.
The Euclidean norm is represented by $\Vert u \Vert \coloneqq \sqrt{\langle u, u \rangle}$.

Let ${\cal M}$ be a submanifold of ${\mathbb{R}^n}$.
A mapping $\gamma \colon (-\varepsilon, \varepsilon) \to {\cal M}$ with some $\varepsilon > 0$ is called a smooth curve around $x \in {\cal M}$ if $\gamma(0) = x$ holds and $\varphi_{\alpha} \circ \gamma$ is smooth for some chart $(U_{\alpha}, \varphi_{\alpha})$ of $\mathcal{M}$ around $x$. Let $\gamma$ be a smooth curve on ${\cal M}$ with $\gamma(0) = x$. We define a linear operator $\dot{\gamma}(0) \colon \mathfrak{F}_{x}({\cal M}) \to \mathbb{R}$ (i.e., the derivative of $\gamma$ at $0 \in \mathbb{R}$) as
\begin{eqnarray*}
\dot{\gamma}(0)[\phi] \coloneqq \left. \frac{d}{dt} (\phi \circ \gamma)(t) \right|_{t=0},
\end{eqnarray*}
where $\phi \in \mathfrak{F}_{x}({\cal M})$, and $\mathfrak{F}_{x}({\cal M})$ is the set of all real-valued $C^{\infty}$ functions defined in some open neighborhood of $x$. The operator $\dot{\gamma}(0)$ is called a tangent vector to $\gamma$ at $x$.
The set of all possible tangent vectors at $x$ is called the tangent space to ${\cal M}$ at $x$ and is denoted by $T_{x} {\cal M}$.
We define the tangent bundle of ${\cal M}$ as $T {\cal M} \coloneqq \bigsqcup_{x \in {\cal M}} T_{x} {\cal M}$, which also becomes a manifold with a dimension twice that of $\mathcal{M}$.
Let $\mathcal{N}$ be a Riemannian manifold and let $\Phi \colon {\cal M} \to \mathcal{N}$ be a differentiable function. The derivative of $\Phi$ at $x \in {\cal M}$, denoted by $\D \Phi(x)$, is a linear operator from $T_{x} {\cal M}$ to $T_{\Phi(x)} \mathcal{N}$ such that 
\begin{eqnarray*}
\D \Phi(x) [\dot{c}(0)] = \left. \frac{d}{dt} (\Phi \circ c)(t) \right|_{t=0},
\end{eqnarray*}
where $c$ is a curve on ${\cal M}$ such that $c(0) = x$.

The chain rule yields that, for any $\dot{\gamma}(0) \in T_{x} {\cal M}$ for some curve $\gamma$ with $\gamma(0) = x$, there exists $\gamma^{\prime}(0) \coloneqq \frac{d}{dt} \gamma(t) |_{t=0} \in {\mathbb{R}^n}$ such that, for all $\phi \in \mathfrak{F}_{x}({\cal M})$ and its smooth extension $\bar{\phi} \colon \mathbb{R}^n \to \mathbb{R}$ (i.e., $\bar{\phi}|_{\mathcal{M}} = \phi$), the following holds:
\begin{equation*}
\dot{\gamma}(0)[\phi] = \left. \frac{d}{dt} (\phi \circ \gamma)(t) \right|_{t=0} = \langle \nabla \bar{\phi}(\gamma(0)), \gamma^{\prime}(0) \rangle = \langle \nabla \bar{\phi}(x), \gamma^{\prime}(0) \rangle.
\end{equation*}
Hence, the linear mapping ${\mathbb{R}^n} \ni \gamma^{\prime}(0) \mapsto \dot{\gamma}(0) \in T_{x} {\cal M}$ is surjective.
On the other hand, if ${\cal M} = {\mathbb{R}^n}$, then the linear mapping ${\mathbb{R}^n} \ni \gamma^{\prime}(0) \mapsto \dot{\gamma}(0) \in T_{x} {\mathbb{R}^n}$ is injective.
Therefore, we identify $T_{x} {\mathbb{R}^n}$ with ${\mathbb{R}^n}$ through this linear bijection (i.e., isomorphism) and denote $T_{x} {\mathbb{R}^n} \cong {\mathbb{R}^n}$.
In a general case where $\mathcal{M}$ is embedded in $\mathbb{R}^n$, we represent the inclusion mapping by $\iota \colon {\cal M} \to {\mathbb{R}^n}$.
Then, we identify $T_{x} {\cal M}$ with a subspace of ${\mathbb{R}^n}$, that is, using the mapping $\D \iota(x) \colon T_{x} {\cal M} \to T_{x} {\mathbb{R}^n} \cong {\mathbb{R}^n}$, we identify $T_x \mathcal{M}$ with the image of $\D \iota(x)$, considering that $T_x \mathcal{M} \cong \D \iota(x)[T_x \mathcal{M}]$.

Let ${\cal M}$ be equipped with a Riemannian metric $x \mapsto \langle \cdot, \cdot \rangle_{x}$ such that for each $x \in {\cal M}$, $\langle \xi, \zeta \rangle_{x} \coloneqq \langle \xi, \zeta \rangle = \xi^\top \zeta$ for all $\xi, \, \zeta \in T_{x} {\cal M} \subset T_{x} {\mathbb{R}^n} \cong {\mathbb{R}^n}$, which is the induced metric from the standard Euclidean metric in $\mathbb{R}^n$.
Then, $\cal M$ is a Riemannian submanifold of $\mathbb{R}^n$.
The norm of $\xi \in T_{x} {\cal M}$ is defined by $\Vert \xi \Vert_{x} \coloneqq \sqrt{ \langle \xi, \xi \rangle_{x} } = \sqrt{\langle \xi, \xi \rangle} = \Vert \xi \Vert$, where $\|\cdot\|$ denotes the $2$-norm (Euclidean norm).
Therefore, the inner products and norms of tangent vectors to $\mathcal{M}$ coincide with those of the vectors regarded as numerical vectors in $\mathbb{R}^n$. 
The norm of a linear operator $A \colon T_{x}{\cal M} \to T_{x}{\cal M}$ is defined by $\Vert A \Vert \coloneqq \sup \{ \Vert A [v] \Vert_{x} \mid \Vert v \Vert_{x} = 1 \}$.

Let $x \in {\cal M}$, and let $\psi \colon {\cal M} \to \mathbb{R}$ be a smooth function on $\mathcal{M}$. The Riemannian gradient of $\psi$ at $x$ is denoted by $\grad \psi(x)$, which is defined as the unique tangent vector at $x$ such that, for any $\xi \in T_x \cal M$,
\begin{equation}
\langle \grad f(x), \xi\rangle_x = \D f(x)[\xi]
\end{equation}
holds.
The Riemannian Hessian of $\psi$ at $x$ is denoted as $\Hess \psi(x)$.\footnote{Here, the Hessian is defined through the Levi-Civita connection on $\cal M$. Because the proposed line search handles only the Riemannian gradient, we have not rigorously defined the Riemannian Hessian herein. Further details are available in the literature~\cite{AbsMahSep2008,boumal2023introduction}.}
Moreover, in the Euclidean space, $\grad \psi(x)$ and $\Hess \psi(x)$ are written as $\nabla \psi(x)$ and $\nabla^{2} \psi(x)$, which are the Euclidean gradient and Hessian, respectively.
Notably, we assume that ${\cal M}$ is an embedded manifold of ${\mathbb{R}^n}$.
Hence, for smooth function $\psi \colon \mathbb{R}^n \to \mathbb{R}$, we denote the Riemannian gradient of the restriction of $\psi$ to $\mathcal{M}$ at $x \in \mathcal{M}$ (i.e., $\grad \psi|_{\mathcal{M}}(x)$) by $\grad \psi(x)$, with a slight abuse of notation.

\subsection{Existing concepts and techniques in Riemannian optimization} \label{sec:pre}

In this section, we consider solving the following unconstrained optimization problem on a Riemannian manifold ${\cal M}$:
\begin{equation}
\begin{array}{ll}
\displaystyle \Mini_{x \in {\cal M}} \ & f(x),
\end{array}
\label{UNPonM}
\end{equation}
where $f \colon {\mathbb{R}^n} \to \mathbb{R}$ is continuously differentiable in ${\mathbb{R}^n}$. Moreover, we provide a prototype algorithm to solve Problem~\eqref{UNPonM}.
To this end, we review some existing concepts and techniques associated with Riemannian optimization.

First, the optimality condition for \eqref{UNPonM} is as follows:
\begin{definition}
We say that $x \in {\cal M}$ satisfies the first-order necessary condition for Problem~\eqref{UNPonM} if $\grad f(x) = 0$. Moreover, we call a point $x \in {\cal M}$ that satisfies $\grad f(x) = 0$ a stationary point.
\end{definition}
Most existing optimization methods for Riemannian optimization problems as well as Euclidean ones are developed to find a stationary point.
The purpose of the prototype algorithm described in this section is also to find a stationary point.

Generally, directly obtaining a stationary point for \eqref{UNPonM} is difficult.
Therefore, the development of an iterative method that generates a sequence converging to a stationary point is considered.
Indeed, several types of iterative methods, such as gradient-type~\cite{sato2015new,sato2019riemannian,zhu2020riemannian} and Newton-type methods~\cite{hu2018adaptive,li2021nonmonotone}, have been proposed for unconstrained optimization problems on Riemannian manifolds.
Most optimization methods can be divided into two types based on their convergence properties: globally convergent and locally fast convergent.
The globally convergent methods have the property that, from an arbitrary initial point, a sequence converging to a stationary point is generated or the norms of generated points approach $0$.
The locally fast convergent methods have the property that a sequence rapidly converging to a stationary point is generated from an initial point sufficiently close to the stationary point.
In this study, we focus on globally convergent methods.
In the case of ${\cal M} = {\mathbb{R}^n}$, such types of methods consist of the following three steps:
\begin{description}
\item[Step 1.] Compute a search direction $p_{k} \in {\mathbb{R}^n}$.
\item[Step 2.] Determine a step size $\alpha_{k} > 0$.
\item[Step 3.] Update a current point $x_{k} \in {\mathbb{R}^n}$ as $x_{k+1} \coloneqq x_{k} + \alpha_{k} p_{k}$.
\end{description}
However, Step 3 may not be defined in the case of a general manifold ${\cal M}$ because the addition operation may not make sense.
To extend the aforementioned framework to Riemannian optimization, the following ideas are typically used: obtaining the search direction $p_{k}$ belonging to a tangent space $T_{x_{k}} {\cal M}$ and computing the next point $x_{k+1} \in {\cal M}$ on the geodesic or its approximation emanating from the current point $x_{k}$ in the direction of $p_{k}$.
To perform these operations, a map that retracts $p_{k}$ onto ${\cal M}$ is required. Therefore, we use the following retraction map introduced in~\cite{AbsMahSep2008}:
\begin{definition}
\label{def:retraction}
We say that $R \colon T{\cal M} \to {\cal M}$ is a retraction on ${\cal M}$ if $R$ is a smooth map and, for any $x \in \mathcal{M}$, the restriction of $R$ to $T_{x} {\cal M}$, denoted by $R_{x}$, satisfies the following properties:
\begin{enumerate}
\item $R_{x}(0_{x}) = x$, where $0_{x}$ denotes the zero vector of $T_{x} {\cal M}$.
\item With the canonical identification $T_{0_{x}} T_{x} {\cal M} \cong T_{x} {\cal M}$, the map $R_{x}$ satisfies
\begin{eqnarray*}
\D R_{x} (0_{x}) = \id_{T_{x} {\cal M}},
\end{eqnarray*}
where $\id_{T_{x} {\cal M}}$ denotes the identity map on $T_{x} {\cal M}$.
\end{enumerate}
\end{definition}
In summary, a prototype algorithm for solving Problem~\eqref{UNPonM} is as follows:
\begin{algorithm}[tbh]
\caption{Prototype algorithm for Problem~\eqref{UNPonM}} \label{prototype_algo}
\begin{algorithmic}[1]
\Require 
Choose an initial point $x_{0} \in {\cal M}$.
\For{$k = 0, 1, 2, \ldots$}
\State{Compute a search direction $p_{k} \in T_{x_{k}} {\cal M}$.} \Comment{Step~1}
\State{Determine a step size $\alpha_{k} > 0$.} \Comment{Step~2} 
\State{Update a current point $x_{k} \in {\cal M}$ as $x_{k+1} \coloneqq R_{x_{k}}(\alpha_{k} p_{k})$.} \Comment{Step~3}
\EndFor
\end{algorithmic}
\end{algorithm}

\section{Armijo line search in Riemannian optimization and its computational cost improvement}
\label{sec:Armijo}
In this section, we provide an overview of the Armijo line search and discuss its challenges related to computational cost.
Moreover, to address these challenges, we propose a class of Riemannian optimization algorithms with an improved Armijo line search to solve Problem~\eqref{UNPonM}. 

\subsection{Overview of the Armijo line search}
Line search plays a crucial role in globally convergent optimization methods.
In the field of optimization in Euclidean spaces, various line search strategies have been proposed.
Some of them, such as the Armijo and Wolfe conditions, have been extended to Riemannian optimization. This section mainly focuses on the Armijo line search in Riemannian optimization.

In this subsection, we provide a brief overview of the line search.
We assume that $x_{k}$ is a current point with $\grad f(x_{k}) \not = 0$ and $p_{k}$ is a search direction satisfying $\langle \grad f(x_{k}), p_{k} \rangle_{x_{k}} < 0$, that is, $p_k$ is a descent direction at $x_k$.
Then, the Armijo line search finds a step size $\alpha_{k}$ satisfying the following Armijo condition:
\begin{eqnarray} \label{Armijo_Riemann}
f(R_{x_{k}}(\alpha_{k} p_{k})) \leq f(x_{k}) + \tau \alpha_{k} \langle \grad f(x_{k}), p_{k} \rangle_{x_{k}},
\end{eqnarray}
where $\tau \in (0,1)$ is a constant.
We note that the existence of $\alpha_{k} > 0$ satisfying \eqref{Armijo_Riemann} can be easily proved~\cite{sato2021riemannian}.
In other words, this method can be used for obtaining a new point $x_{k+1} \coloneqq R_{x_{k}}(\alpha_{k} p_{k}) \in {\cal M}$ such that $f(x_{k+1}) < f(x_{k})$.
Indeed, we have
\begin{equation}
f(x_{k+1}) = f(R_{x_{k}}(\alpha_{k} p_{k})) \leq f(x_{k}) + \tau \alpha_{k} \langle \grad f(x_{k}), p_{k} \rangle_{x_{k}} < f(x_{k})
\end{equation}
because of $\langle \grad f(x_{k}), p_{k} \rangle_{x_{k}} < 0$. Moreover, if ${\cal M} = \mathbb{R}^{n}$, then condition~\eqref{Armijo_Riemann} is reformulated as
\begin{eqnarray} \label{Armijo_Euclidean}
f(x_{k} + \alpha_{k} p_{k}) \leq f(x_{k}) + \tau \alpha_{k} \langle \nabla f(x_{k}), p_{k} \rangle.
\end{eqnarray}
When implementing the Armijo line search on a computer, the following algorithm (Algorithm~\ref{Armijo_algo}) based on backtracking is typically used.
\begin{algorithm}[htbp]
\caption{Armijo line search algorithm based on backtracking} \label{Armijo_algo}
\begin{algorithmic}[1]
\Require Set constants $\beta, \, \tau \in (0,1)$, and $\bar{\alpha} > 0$.
\Function{Armijo-line search}{$x_{k}, p_{k}$}
\For{$l = 0, 1, 2, \ldots$} \label{line2}
\If{$f(R_{x_{k}}(\beta^{l} \bar{\alpha}p_{k})) \leq f(x_{k}) + \tau \beta^{l}\bar{\alpha} \langle \grad f(x_{k}), p_{k} \rangle_{x_{k}}$}
\State{Set $l_{k} \coloneqq l$ and break.}
\EndIf
\EndFor \label{line4}
\State{\Return $\alpha_{k} \coloneqq \beta^{l_{k}}\bar{\alpha}$.}
\EndFunction
\end{algorithmic}
\end{algorithm}

For a Riemannian manifold $\cal M$, \eqref{Armijo_Riemann} is a natural condition that leads to a reasonable reduction of the objective function value.
However, when $\cal M$ is a Riemannian submanifold of $\mathbb{R}^n$, condition \eqref{Armijo_Euclidean} also makes sense, although it is important to investigate if~\eqref{Armijo_Euclidean} actually holds in reducing the function value; this is discussed in the next subsection.
Here, we consider the differences in computational cost between \eqref{Armijo_Riemann} and \eqref{Armijo_Euclidean}, that is, the Riemannian and Euclidean Armijo conditions.
Generally, the computational cost of the line search for the Riemannian Armijo condition is higher than that for the Euclidean one because the left-hand side of \eqref{Armijo_Riemann} requires the calculation of the retraction.
If $l_k$ is large, that is, the loop-statement of Lines~\ref{line2}--\ref{line4} in Algorithm~\ref{Armijo_algo} is performed many times and the computational cost regarding the retraction is large, then the differences in their computational costs increase.
This phenomenon is a notable drawback of the Armijo line search in Riemannian optimization.
The same phenomenon also applies to other line search strategies, for example, the Wolfe line search in Riemannian optimization.

\subsection{Modification of the Armijo line search in Riemannian optimization}
Hereafter, we refer to the line search based on~\eqref{Armijo_Riemann} as the Riemannian Armijo line search and~\eqref{Armijo_Euclidean} as the Euclidean Armijo line search.
As mentioned above, the Riemannian Armijo line search has a drawback in terms of the computational cost associated with the retraction.
This subsection discusses how this study addresses the aforementioned drawback.

\begin{figure}
\centering
\includegraphics[width=0.45\linewidth]{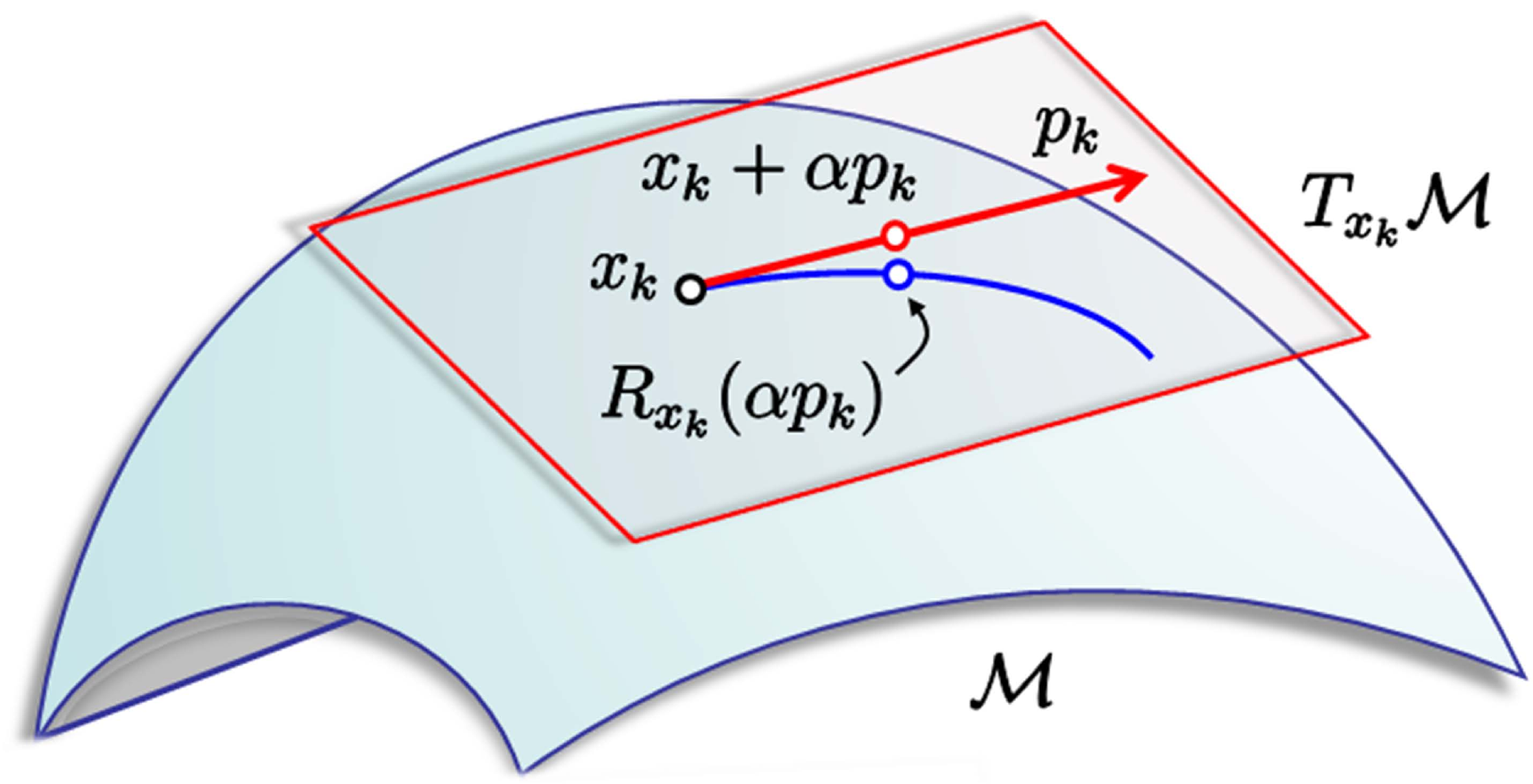}
    \caption{Concept of the proposed line search strategy. We approximate $R_{x_k}(\alpha p_k)$ by $x_k + \alpha p_k$.}
    \label{fig:concept}
\end{figure}

The idea of this study is as follows.
A Riemannian submanifold of the Euclidean space is generally a curved surface.
Therefore, the Euclidean Armijo line search cannot be applied directly.
However, around each point, the (smooth) manifold can be locally and linearly approximated by its tangent space at the point.
Since the search direction is computed as a tangent vector at the current point, the Euclidean Armijo line search along this direction means the line search on the tangent space. Therefore, the Euclidean Armijo line search has the potential to well approximate the Riemannian Armijo line search. The concept is illustrated in Figure~\ref{fig:concept}.

Before considering the specific improvement strategy of the Riemannian Armijo line search, we highlight an important property of the gradient of $f$ that relates the Riemannian Armijo condition to the Euclidean one.
To this end, let ${\cal P}_{x} \colon {\mathbb{R}^n} \to T_{x} {\cal M}$ be the orthogonal projection; we prepare the following lemma that provides a relationship between $\grad f(x)$ and $\nabla f(x)$.
The lemma is a well-known result~\cite{AbsMahSep2008}; however, for self-containedness and a better understanding of the proposed strategy, we provide a proof as follows.
\begin{lemma} \label{lemma:grad}
Let $\cal M$ be a Riemannian submanifold of $\mathbb{R}^n$ and let $x \in {\cal M}$ and $p \in T_{x} {\cal M}$.
Then, $\langle \grad f(x), p \rangle_{x} = \langle \nabla f(x), p \rangle$ and $\grad f(x) = {\cal P}_{x}(\nabla f(x))$ hold.
\end{lemma}

\begin{proof}
Let $p \in T_{x} {\cal M}$ be an arbitrary tangent vector.
Then, there exists a curve $c \colon (-\varepsilon, \varepsilon) \to {\cal M}$ such that $c(0) = x$ and $c^{\prime}(0) = p$. Exploiting the relationship ${\cal M} \subset {\mathbb{R}^n}$ and chain rule yields
\begin{eqnarray*}
\langle \grad f(x), p \rangle_{x} &=& \D f(x)[p] = \left. \frac{d(f \circ c)(t)}{dt} \right|_{t=0} = \langle \nabla f(c(0)), c^{\prime}(0) \rangle = \langle \nabla f(x), p \rangle.
\end{eqnarray*}
It follows that $\langle \grad f(x), p \rangle_{x} = \langle \nabla f(x), p \rangle$. Moreover, by considering the orthogonal decomposition of $\nabla f(x)$, we have
\begin{eqnarray*}
\langle \grad f(x), p \rangle_{x} &=& \langle \nabla f(x), p \rangle
\\
&=& \langle {\cal P}_{x}(\nabla f(x)), p \rangle + \langle {\cal P}_{x}^{\bot}(\nabla f(x)), p \rangle 
\\
&=& \langle {\cal P}_{x}(\nabla f(x)), p \rangle_{x},
\end{eqnarray*}
where ${\cal P}_{x}^{\bot} \colon {\mathbb{R}^n} \to T_{x} {\cal M}^{\bot}$ denotes the orthogonal projection onto the normal space $T_x \mathcal{M}^\perp$, which is the orthogonal complement of $T_x \mathcal{M}$ in $\mathbb{R}^n$ with respect to the inner product $\langle \cdot, \cdot\rangle$.
Consequently, we obtain $\langle \grad f(x), p \rangle_{x} = \langle {\cal P}_{x}(\nabla f(x)), p \rangle_{x}$ for all $p \in T_{x} {\cal M}$, that is, $\grad f(x) = {\cal P}_{x}(\nabla f(x))$.
Therefore, the claim is proved. 
\end{proof}
As evident from Lemma~\ref{lemma:grad}, the right-hand sides of the Riemannian Armijo condition~\eqref{Armijo_Riemann} and the Euclidean one~\eqref{Armijo_Euclidean} coincide.

Furthermore, using Lemma~\ref{lemma:grad}, we provide a property of $f$ regarding the left-hand sides of~\eqref{Armijo_Riemann} and~\eqref{Armijo_Euclidean}, which plays a crucial role in improving the Riemannian Armijo line search.
\begin{proposition} \label{prop:ErrorResidue}
Let $\cal M$ be a Riemannian submanifold of $\mathbb{R}^n$ and let $x \in {\cal M}$ and $p \in T_{x} {\cal M}$ be given. Moreover, let the function $E \colon \mathbb{R} \to \mathbb{R}$ be defined by $E(\alpha) \coloneqq | f(x + \alpha p) - f(R_{x}(\alpha p)) |$.
Then, $E(\alpha) = o(\alpha) ~ (\alpha \to 0)$ is satisfied.
\end{proposition}
\begin{proof}
From the first-order Taylor expansions of $f(x + \alpha p)$ and $f(R_{x}(\alpha p))$, we have the following:
\begin{eqnarray*}
f(x + \alpha p) &=& f(x) + \alpha \langle \nabla f(x), p \rangle + \alpha r_{1}(\alpha),
\\
f(R_{x}(\alpha p)) &=& f(x) + \alpha \langle \grad f(x), p \rangle_{x} + \alpha r_{2}(\alpha),
\end{eqnarray*}
where $r_{1}$ and $r_{2}$ are some functions such that $r_{1}(\alpha) \to 0$ and $r_{2}(\alpha) \to 0$ as $\alpha \to 0$.
From these equalities and Lemma~\ref{lemma:grad}, the assertion is proved as follows:
\begin{equation}
\bigg|\frac{E(\alpha)}{\alpha}\bigg| = |r_1(\alpha) - r_2(\alpha)| \to 0 \quad (\alpha \to 0).
\end{equation}
\end{proof}
Proposition~\ref{prop:ErrorResidue} states that $f(x + \alpha p)$ appropriately approximates $f(R_{x}(\alpha p))$ when $\alpha > 0$ is sufficiently small.
Thus, we expect that a globally convergent method for Problem~\eqref{UNPonM} can be developed even if another line search based on the Euclidean, instead of Riemannian, Armijo condition~\eqref{Armijo_Euclidean} is used.

Therefore, as Algorithm~\ref{alg:Newton_method}, we propose a class of optimization methods for Problem~\eqref{UNPonM} equipped with a new line search based on \eqref{Armijo_Euclidean} to overcome the aforementioned drawback of the Riemannian Armijo line search.
We also refer to this Euclidean Armijo line search as the modified (Riemannian) Armijo line search.
\begin{algorithm}[tbh]
\caption{Class of optimization methods for Problem~\eqref{UNPonM} with modified Riemannian Armijo line search} \label{alg:Newton_method}
\begin{algorithmic}[1]
\Require 
Set constants $\beta, \, \tau \in (0,1)$, and $\bar{\alpha} > 0$. Choose an initial point $x_{0} \in {\cal M}$.
\For{$k = 0, 1, 2, \ldots$}
\State{Compute a search direction $p_{k} \in T_{x_{k}} {\cal M}$ satisfying} \Comment{Step~1}
\begin{eqnarray*}
H_{k} [p_{k}] = - \grad f(x_{k}).
\end{eqnarray*}
\For{$l = 0, 1, 2, \ldots$} \Comment{Step~2}
\If{$f(x_{k} + \beta^{l}\bar{\alpha} p_{k}) \leq f(x_{k}) + \tau \beta^{l}\bar{\alpha} \langle \grad f(x_{k}), p_{k} \rangle_{x_{k}}$} \label{Line:FirstIf}
\If{$f(R_{x_{k}}(\beta^{l}\bar{\alpha} p_{k})) \leq f(x_{k}) + \tau \beta^{l}\bar{\alpha} \langle \grad f(x_{k}), p_{k} \rangle_{x_{k}}$} \label{Line:SecondIf}
\State{Set $l_{k} \coloneqq l$ and break.}
\EndIf
\EndIf
\EndFor
\State{Determine a step size as $\alpha_{k} \coloneqq \beta^{l_{k}}\bar{\alpha}$.}
\State{Update a current point $x_{k} \in {\cal M}$ as $x_{k+1} \coloneqq R_{x_{k}}(\alpha_{k} p_{k})$.} \Comment{Step~3}
\EndFor
\end{algorithmic}
\end{algorithm}
In the algorithm, the search direction $p_k \in T_{x_k} \cal M$ is computed as the solution to the linear equation
\begin{equation}
H_k[p_k] = -\grad f(x_k).
\end{equation}
Therefore, theoretically, we have $p_k = -H_k^{-1} \grad f(x_k)$, where $H_{k} \colon T_{x_{k}} {\cal M} \to T_{x_{k}} {\cal M}$ is an invertible linear operator in $T_{x_k} \cal M$ with self-adjointness, that is, 
\begin{eqnarray*}
\langle H_{k}[u], v \rangle_{x_{k}} = \langle u, H_{k}[v] \rangle_{x_{k}} \text{ for all } u, v \in T_{x_{k}} {\cal M}.
\end{eqnarray*}
This choice of search direction covers a wide class of Riemannian optimization algorithms because $H_k$ involves a certain degree of arbitrariness.
In other words, a specific choice of $H_k$ corresponds to a specific optimization method, for example,
\begin{itemize}
\item
$H_k = \id_{T_{x_k}\cal M}$ $\implies$ steepest descent direction $p_k = -\grad f(x_k)$;
\item 
$H_k = \Hess f(x_k)$ $\implies$ Newton's direction, which is the solution of the Newton equation $\Hess f(x_k)[p_k] = -\grad f(x_k)$;
\item 
$H_k$ satisfies the secant condition $\implies$ quasi-Newton's direction.
\end{itemize}

Algorithm~\ref{alg:Newton_method} adopts a new line search strategy based on the Euclidean Armijo condition~\eqref{Armijo_Euclidean} to determine the step size.
In Step~2 of Algorithm~\ref{alg:Newton_method}, backtracking is performed with the initial guess of the step length, $\bar{\alpha} > 0$, to find a step length satisfying the Euclidean (not Riemannian) Armijo condition, with the aim of reasonable reduction of computational cost by avoiding retractions.
After finding an integer $l$ such that the candidate for the step length $\beta^l\bar{\alpha}$ satisfies the Euclidean Armijo condition (noting that $\langle \nabla f(x_k), p_k\rangle = \langle \grad f(x_k), p_k\rangle_{x_k}$) as
\begin{equation}
\label{eq:alg_Euc_Armijo}
f(x_k + \beta^l\bar{\alpha} p_k) \leq f(x_k) + \tau \beta^l\bar{\alpha} \langle \grad f(x_k), p_k\rangle_{x_k},
\end{equation}
we check if this $\beta^l \bar{\alpha}$ satisfies the Riemannian Armijo condition as
\begin{equation}
\label{eq:alg_R_Armijo}
f(R_{x_k}(\beta^l \bar{\alpha} p_k)) \leq f(x_k) + \tau \beta^l \bar{\alpha} \langle \grad f(x_k), p_k\rangle_{x_k}
\end{equation}
because we finally need a step length satisfying the Riemannian Armijo condition.
However, if the candidate $\beta^l \bar{\alpha}$ satisfying~\eqref{eq:alg_Euc_Armijo} is small, that is, if $l$ is large, we can expect that the candidate $\beta^l \bar{\alpha}$ tends to satisfy~\eqref{eq:alg_R_Armijo} since the Euclidean Armijo condition~\eqref{eq:alg_Euc_Armijo} and the Riemannian one~\eqref{eq:alg_R_Armijo} are close to each other as discussed.

In the next section, we prove the global convergence of Algorithm~\ref{alg:Newton_method} under some appropriate assumptions.

\begin{remark}
The initial guess of the step length, $\bar{\alpha} > 0$, can theoretically be any value.
A natural practical choice is $\bar{\alpha} = 1$.
In particular, for the Newton method ($H_k = \Hess f(x_k)$), $\bar{\alpha} = 1$ should first be attempted to achieve local fast convergence.
\end{remark}

\section{Global convergence of Algorithm~\ref{alg:Newton_method}}
\label{sec:Convergence}
We make certain assumptions to guarantee the global convergence of Algorithm~\ref{alg:Newton_method}.
\begin{assumption} \label{assumption:global}
~
\begin{description}
\item[{\rm (A1)}] 
The $C^1$ function $\bar{f} \colon \mathbb{R}^n \to \mathbb{R}$ is an extension of the objective function $f \colon \mathcal{M} \to \mathbb{R}$, and we also denote $\bar{f}$ as $f$ with a slight abuse of notation.
The (Euclidean) gradient of the function $f \colon \mathbb{R}^n \to \mathbb{R}$ is Lipschitz continuous, that is,
\begin{eqnarray*}
\Vert \nabla f(x) - \nabla f(y) \Vert \leq L \Vert x - y \Vert \text{ for all } x, y \in {\mathbb{R}^n},
\end{eqnarray*}
where $L > 0$ is a constant.
\item[{\rm (A2)}] A sequence $\{ x_{k} \}$ generated by Algorithm~{\rm \ref{alg:Newton_method}} is bounded.
\item[{\rm (A3)}] There exist positive constants $\nu$ and $\rho$ such that, for each $k \in \mathbb{N} \cup \{ 0 \}$,
\begin{eqnarray*}
\nu \Vert u \Vert_{x_{k}}^{2} \leq \langle H_{k}[u], u \rangle_{x_{k}} \leq \rho \Vert u \Vert_{x_{k}}^{2} \text{ for all } u \in T_{x_{k}} {\cal M},
\end{eqnarray*}
where $\{ x_{k} \}$ is a sequence generated by Algorithm~{\rm \ref{alg:Newton_method}}.
\end{description}
\end{assumption}
If $\grad f(x_k) = 0$ is satisfied for some $k \in \mathbb{N} \cup \{0\}$, then the goal is considered to be achieved.
Therefore, in the following argument, we assume that Algorithm~\ref{alg:Newton_method} generates an infinite sequence $\{ x_{k} \}$ that satisfies $\grad f(x_{k}) \not = 0$ for all $k \in \mathbb{N} \cup \{ 0 \}$.

This section is divided into two parts.
In the first part, we discuss the well-definedness of Algorithm~\ref{alg:Newton_method}, and in the second part, we prove its global convergence. 

\subsection{Well-definedness of Algorithm~\ref{alg:Newton_method}}
First, we verify that, given $x_k \in \cal M$ and $p_k \in T_{x_k} \cal M$, Step~1 of Algorithm~\ref{alg:Newton_method} can find the unique solution of the equation $H_k[p_k] = -\grad f(x_k)$ for $p_k \in T_{x_k}\mathcal{M}$.
Here, we note that $T_{x_{k}} {\cal M}$ is a finite-dimensional subspace of the Euclidean space ${\mathbb{R}^n}$, and $H_{k} \colon T_{x_{k}} {\cal M} \to T_{x_{k}} {\cal M}$ is a linear operator with self-adjointness.
Since (A3) of Assumption~\ref{assumption:global} holds, the operator $H_{k}$ is positive definite and, therefore, invertible.
Hence, the equation has the unique solution $p_{k} = -H_{k}^{-1} \grad f(x_{k})$.

To show the well-definedness of Algorithm~\ref{alg:Newton_method}, we further need to ensure the existence of the integer $l_{k}$ calculated in Step~2 (i.e., the finite termination of Step~2).
For this purpose, we provide the following proposition.
\begin{proposition} \label{well-def:Armijo}
Suppose that {\rm (A1)} and {\rm (A3)} of Assumption~{\rm \ref{assumption:global}} hold.
Let $\delta \coloneqq 2\nu(1-\tau)/L$.
For any $\alpha \in (0, \delta]$, the inequality 
\begin{equation}
f(x_{k} + \alpha p_{k}) \leq f(x_{k}) + \tau \alpha \langle \grad f(x_{k}), p_{k} \rangle_{x_{k}}
\end{equation}
holds.
\end{proposition}
\begin{proof}
Let $\alpha \in (0, \delta]$.
Since (A1) of Assumption~\ref{assumption:global} holds, it follows from \cite[Proposition~A.24]{Ber99} that
\begin{eqnarray}
f(x_{k} + \alpha p_{k}) \leq f(x_{k}) + \alpha \langle \nabla f(x_{k}), p_{k} \rangle + \frac{L \alpha^{2}}{2} \Vert p_{k} \Vert^{2}. \label{ineq:Lip}
\end{eqnarray}
Further, $\langle \nabla f(x_{k}), p_{k} \rangle$ in \eqref{ineq:Lip} can be evaluated as follows:
\begin{align}
\langle \nabla f(x_{k}), p_{k} \rangle
&= \langle \grad f(x_k), p_k\rangle_{x_k}\\
&=\tau \langle \grad f(x_{k}), p_{k} \rangle_{x_{k}} + (1 - \tau) \langle \grad f(x_{k}), p_{k} \rangle_{x_{k}} \nonumber
\\
&= \tau \langle \grad f(x_{k}), p_{k} \rangle_{x_{k}} - (1 - \tau) \langle H_{k} p_{k}, p_{k} \rangle_{x_{k}} \nonumber
\\
&\leq \tau \langle \grad f(x_{k}), p_{k} \rangle_{x_{k}} - \nu (1 - \tau) \Vert p_{k} \Vert_{x_{k}}^{2}, \label{ineq:Armijo}
\end{align} 
where the first and third equalities follow from Lemma~\ref{lemma:grad} and the equation $H_{k} [p_{k}] = -\grad f(x_{k})$, respectively, and the inequality is derived from (A3) of Assumption~\ref{assumption:global}. By combining \eqref{ineq:Lip} and \eqref{ineq:Armijo}, we derive
\begin{eqnarray*}
f(x_{k} + \alpha p_{k}) 
&\leq& f(x_{k}) + \tau \alpha \langle \grad f(x_{k}), p_{k} \rangle_{x_{k}} + \alpha \Vert p_{k} \Vert_{x_{k}}^{2} \left\{ \frac{L \alpha}{2} - \nu (1 - \tau) \right\}
\\
&\leq& f(x_{k}) + \tau \alpha \langle \grad f(x_{k}), p_{k} \rangle_{x_{k}},
\end{eqnarray*}
because of $L\alpha/2 - \nu(1-\tau) \leq 0$ from $0 < \alpha \leq \delta = 2\nu(1-\tau)/L$.
Therefore, the desired inequality is obtained.
\end{proof}

From Proposition~\ref{well-def:Armijo}, the condition described in Line~\ref{Line:FirstIf} of Algorithm~\ref{alg:Newton_method} is satisfied for a sufficiently large $l \in \mathbb{N} \cup \{ 0 \}$ such that $\beta^{l} \in (0, \delta]$. Moreover, the condition described in Line~\ref{Line:SecondIf} is the Riemannian Armijo line search; therefore, it is clear that a sufficiently large $l \in \mathbb{N} \cup \{ 0 \}$ satisfies this condition. Hence, Algorithm~\ref{alg:Newton_method} is well-defined under Assumption~\ref{assumption:global}.

\subsection{Global convergence of Algorithm~\ref{alg:Newton_method}}
Here, we show the global convergence of Algorithm~\ref{alg:Newton_method}.
We first provide a lemma for the sequence of search directions. 

\begin{lemma} \label{lemma:pbounded}
Suppose that Assumption~{\rm \ref{assumption:global}} holds. Then, the sequence $\{ p_{k} \}$, which is generated by Step~{\rm 1} of Algorithm~{\rm \ref{alg:Newton_method}}, is bounded.
\end{lemma}
\begin{proof}
From (A3) of Assumption~\ref{assumption:global}, the equation $H_{k} [p_{k}] = -\grad f(x_{k})$, and the Cauchy--Schwartz inequality, we obtain
\begin{eqnarray*}
\nu \Vert p_{k} \Vert_{x_{k}}^{2} \leq \langle H_{k} [p_{k}], p_{k} \rangle_{x_{k}} = -\langle \grad f(x_{k}), p_{k} \rangle_{x_{k}} \leq \Vert \grad f(x_{k}) \Vert_{x_{k}} \Vert p_{k} \Vert_{x_{k}}
\end{eqnarray*}
for all $k \in \mathbb{N} \cup \{ 0 \}$.
The boundedness of $\{ x_{k} \}$ and the continuity of $\grad f$ ensure that there exists a constant $C > 0$ such that $\Vert \grad f(x_{k}) \Vert_{x_{k}} \leq C$ for all $k \in \mathbb{N} \cup \{ 0 \}$.
Therefore, these results imply that $\Vert p_{k} \Vert_{x_{k}} \leq C / \nu$ for all $k \in \mathbb{N} \cup \{ 0 \}$, that is, $\{ p_{k} \}$ is bounded.
\end{proof}
Using this lemma, we prove the following main theorem.
\begin{theorem}
Suppose that Assumption~{\rm \ref{assumption:global}} holds. If Algorithm~{\rm \ref{alg:Newton_method}} generates an infinite sequence $\{ x_{k} \}$ such that $\grad f(x_{k}) \not = 0$ for all $k \in \mathbb{N} \cup \{ 0 \}$, then $\liminf_{k \to \infty} \Vert \grad f(x_{k}) \Vert_{x_{k}} = 0$.
\end{theorem}

\begin{proof}
From Steps~1--3 of Algorithm~\ref{alg:Newton_method}, we observe $H_{k} [p_{k}] = -\grad f(x_{k})$ and $f(x_{k+1}) \leq f(x_{k}) + \tau \beta^{l_{k}}\bar{\alpha} \langle \grad f(x_{k}), p_{k} \rangle_{x_{k}}$ for each $k \in \mathbb{N} \cup \{ 0 \}$, leading to $\tau \beta^{l_{k}}\bar{\alpha} \langle H_{k}[p_{k}], p_{k} \rangle_{x_{k}} = -\tau \beta^{l_{k}}\bar{\alpha} \langle \grad f(x_{k}), p_{k} \rangle_{x_{k}} \leq f(x_{k}) - f(x_{k+1})$.
Then, it follows from (A3) of Assumption~\ref{assumption:global} that
\begin{eqnarray*}
0 \leq \nu \tau \beta^{l_{k}}\bar{\alpha} \Vert p_{k} \Vert_{x_{k}}^{2} \leq f(x_{k}) - f(x_{k+1})
\end{eqnarray*}
for all $k \in \mathbb{N} \cup \{ 0 \}$.
Since the sequence $\{ f(x_{k}) \}$ is bounded below and monotonically decreasing, it converges, implying $\lim_{k \to \infty} (f(x_k) - f(x_{k+1})) = 0$.
Therefore, we have $\beta^{l_{k}} \Vert p_{k} \Vert_{x_{k}}^{2} \to 0$ as $k \to \infty$.
In the following, we prove $\liminf_{k \to \infty} \Vert p_{k} \Vert_{x_{k}} = 0$.

Two possibilities exist: $\liminf_{k \to \infty} \beta^{l_{k}} > 0$ and $\liminf_{k \to \infty} \beta^{l_{k}} = 0$.
In the first case, we readily obtain $\Vert p_{k} \Vert_{x_{k}} \to 0$ as $k \to \infty$.
We proceed to the second case.
We can take ${\cal L} \subset \mathbb{N}$ such that $l_{k} \to \infty$ as ${\cal L} \ni k \to \infty$. From~(A2) of Assumption~\ref{assumption:global} and Lemma~\ref{lemma:pbounded}, there exist $x^{\ast} \in {\mathbb{R}^n}$, $p^{\ast} \in {\mathbb{R}^n}$, and ${\cal K} \subset {\cal L}$ such that $x_{k} \to x^{\ast}$ and $p_{k} \to p^{\ast}$ as ${\cal K} \ni k \to \infty$.
Since $l_{k} \to \infty$ as ${\cal K} \ni k \to \infty$, there exists $n_{0} \in \mathbb{N}$ such that $l_{k} > \log_{\beta} (\delta / \bar{\alpha}) + 1$ for all $k \in {\cal K}$ with $k \geq n_{0}$, where we define $\delta \coloneqq 2 \nu (1-\tau) / L$ as in Proposition~\ref{well-def:Armijo}.
Here, we arbitrarily take $k \in {\cal K}$ with $k \geq n_{0}$.
Note that $\beta \in (0,1)$ and $\log_{\beta}( \beta^{l_{k} - 1}\bar{\alpha}) = l_{k} - 1 + \log_{\beta}\bar{\alpha} > \log_{\beta} \delta$, that is, $\gamma_{k} \coloneqq \beta^{l_{k} - 1}\bar{\alpha} \in (0, \delta]$.
Hence, Proposition~\ref{well-def:Armijo} ensures that
\begin{equation}
f(x_{k} + \gamma_{k} p_{k}) \leq f(x_{k}) + \tau \gamma_{k} \langle \grad f(x_{k}), p_{k} \rangle_{x_{k}},
\end{equation}
namely, the condition in Line~\ref{Line:FirstIf} of Algorithm~\ref{alg:Newton_method} holds.
This implies that the condition in Line~\ref{Line:SecondIf} does not hold, that is,
\begin{eqnarray}
f(x_{k}) + \tau \gamma_{k} \langle \grad f(x_{k}), p_{k} \rangle_{x_{k}} < f(R_{x_{k}}(\gamma_{k} p_{k})) \label{ineq:notArmijo}
\end{eqnarray}
because if Line~5 held, then the $l_{k}$th loop of Step~2 would not be performed.
Let $\varphi_{k}(t) \coloneqq f(R_{x_{k}}(tp_{k}))$ for $t \in \mathbb{R}$.
Then, we can rewrite~\eqref{ineq:notArmijo} as follows:
\begin{eqnarray}
\hspace{-5mm} (\tau - 1) \langle \grad f(x_{k}), p_{k} \rangle_{x_{k}} < \frac{\varphi_{k}(\gamma_{k}) - \varphi_{k}(0)}{\gamma_{k}} - \langle \grad f(x_{k}), p_{k} \rangle_{x_{k}}. \label{ineq:notArmijo2}
\end{eqnarray}
From (A3) of Assumption~\ref{assumption:global} and the equation $H_{k} [p_{k}] = -\grad f(x_{k})$, we can derive
\begin{equation}
\nu (1 - \tau) \Vert p_{k} \Vert_{x_{k}}^{2} \leq (1-\tau)\langle p_k, H_k[p_k]\rangle_{x_k}
= (\tau - 1) \langle \grad f(x_{k}), p_{k} \rangle_{x_{k}}.
\end{equation}
This inequality together with \eqref{ineq:notArmijo2} yields
\begin{eqnarray}
0 \leq \nu (1 - \tau) \Vert p_{k} \Vert_{x_{k}}^{2} < \frac{\varphi_{k}(\gamma_{k}) - \varphi_{k}(0)}{\gamma_{k}} - \langle \grad f(x_{k}), p_{k} \rangle_{x_{k}}, \label{ineq:notArmijo4}
\end{eqnarray}
where we used $0 < \tau < 1$.
It follows from the mean value theorem regarding $\varphi_{k}$ and chain rule that there exists $\theta_k \in (0, 1)$ such that
\begin{align}
\frac{\varphi_{k}(\gamma_{k}) - \varphi_{k}(0)}{\gamma_{k}} &= \varphi'(\theta_k\gamma_k)\\
&= \D f(R_{x_k}(\theta_k\gamma_kp_k))[\D R_{x_k}(\theta_k\gamma_k p_k)[p_k]]\\
&= \langle \grad f(R_{x_k}(\theta_k\gamma_kp_k)), \D R_{x_k}(\theta_k\gamma_kp_k)[p_k]\rangle_{R_{x_k}(\theta_k\gamma_kp_k)}.
\end{align}
Therefore, we have
\begin{align}
0 &\leq \nu (1 - \tau) \Vert p_{k} \Vert_{x_{k}}^{2} \nonumber
\\
&\leq \langle \grad f(R_{x_{k}}(\theta_{k} \gamma_{k} p_{k})), \D R_{x_{k}}(\theta_{k} \gamma_{k} p_{k})[p_{k}] \rangle_{R_{x_k}(\theta_k\gamma_kp_k)} - \langle \grad f(x_{k}), p_{k} \rangle_{x_{k}}.
\label{ineq:notArmijo5}
\end{align}
Recall that $R_{x^{\ast}}(0) = x^{\ast}$, $\D R_{x^{\ast}}(0) = \id_{T_{x^{\ast}} {\cal M}}$, and $\gamma_{k} = \beta^{l_{k}-1}\bar{\alpha} \to 0$ as ${\cal K} \ni k \to \infty$.
Taking the limit in~\eqref{ineq:notArmijo5} over ${\cal K}$, we obtain $0 \leq \nu (1-\tau) \Vert p^{\ast} \Vert_{x^{\ast}}^{2} \leq 0$, that is, $p^{\ast} = 0$, and therefore, $\Vert p_{k} \Vert_{x_{k}} \to 0$ as ${\cal K} \ni k \to \infty$.

Thus, we have verified that $\liminf_{k \to \infty} \Vert p_{k} \Vert_{x_{k}} = 0$ always holds.
On the other hand, (A3) of Assumption~\ref{assumption:global} implies that $\langle H_{k}[u], u \rangle_{x_{k}} \leq \rho$ for all $u \in T_{x_{k}} {\cal M}$ with $\Vert u \Vert_{x_{k}} = 1$; therefore, taking $u$ as a unit eigenvector with respect to the maximum eigenvalue of $H_{k}$, which is denoted by $\lambda_{\max}(H_{k})$, we have $\Vert H_{k} \Vert = \lambda_{\max}(H_{k}) \leq \rho$.
It then follows from the equation $H_{k} [p_{k}] = -\grad f(x_{k})$ that 
\begin{equation}
\Vert \grad f(x_{k}) \Vert_{x_{k}} = \|H_k[p_k]\|_{x_k} \leq \Vert H_{k} \Vert \Vert p_{k} \Vert_{x_{k}} \leq \rho \Vert p_{k} \Vert_{x_{k}}
\end{equation}
holds.
From this inequality and $\liminf_{k \to \infty} \Vert p_{k} \Vert_{x_{k}} = 0$, we conclude that $\liminf_{k \to \infty} \Vert \grad f(x_{k}) \Vert_{x_{k}} = 0$.
\end{proof}

\section{Numerical experiments}
\label{sec:experiements}
In this section, we present numerical results that demonstrate the superiority of the proposed optimization algorithms equipped with the modified Armijo line search (Algorithm~\ref{alg:Newton_method}) over the existing methods with the conventional Riemannian Armijo search (Algorithm~\ref{Armijo_algo}).
The following numerical experiments were performed in double-precision floating-point arithmetic on a computer (Apple M1 Max, 64 GB RAM) equipped with MATLAB R2024a.

Throughout the experiments, the parameter $\tau \in (0, 1)$ in the Armijo condition was set to $10^{-4}$, which is a standard choice~\cite{nocedal2006numerical}.
The contraction rate $\beta \in (0, 1)$ for the backtracking procedure in Algorithms~\ref{Armijo_algo} and~\ref{alg:Newton_method} was set to $0.5$.
Furthermore, the initial guess of the step length was set as $\bar{\alpha} = 1$.

\subsection{Optimization problems on several manifolds and retractions on them}
In this subsection, we introduce three problems on different manifolds used in the subsequent numerical experiments.
Because the proposed method aims to reduce the number of computations of retractions and we would like to observe how this reduction decreases the computational time for optimization, we also discuss the retractions on the considered manifolds.

Let $n, r \in \mathbb{N}$ be natural numbers with $1 \leq r \leq n$.
In the following numerical experiments, we deal with optimization problems on the sphere
\begin{equation}
S^{n-1} = \{x \in \mathbb{R}^n \mid x^\top x = 1\},
\end{equation}
Stiefel manifold
\begin{equation}
\St(r,n) = \{X \in \mathbb{R}^{n \times r} \mid X^\top X = I_r\},
\end{equation}
and the manifold of all $n \times n$ symmetric positive definite matrices (SPD manifold)
\begin{equation}
\Sym_{++}(n) = \{X \in \Sym(n) \mid X \succ 0\},
\end{equation}
where $\Sym(n) = \{X \in \mathbb{R}^{n \times n} \mid X^\top = X\}$, and 
$X \succ 0$ indicates that $X$ is positive definite, that is, $d^\top X d > 0$ for all nonzero $d \in \mathbb{R}^n$.
The dimensions of the manifolds $S^{n-1}$, $\St(r, n)$, and $\Sym_{++}(n)$ are $n-1$, $r(r+1)/2 + r(n-r)$, and $n(n+1)/2$, respectively.

On the sphere $S^{n-1}$, a commonly used retraction is based on the normalization with respect to the $2$-norm, which is defined for $x \in S^{n-1}$ and $p \in T_x S^{n-1} = \{d \in \mathbb{R}^n \mid x^\top d = 1\}$ as
\begin{equation}
R_x(p) = \frac{x+p}{\|x+p\|},
\end{equation}
where $\|\cdot\|$ is the $2$-norm (Euclidean norm) in $\mathbb{R}^n$.
Since a point $x$ on the sphere $S^{n-1}$ and tangent vector $p \in T_x S^{n-1}$ at $x$ are both represented as $n$-dimensional numerical vectors, the $2$-norm of $x + p$ makes sense.
In addition to just computing $x + \alpha p$, by definition of the retraction, computing $R_x(\alpha p) = (x+\alpha p) / \|x+\alpha p\|$ requires extra costs for computing the norm of $x + \alpha p$ and dividing each element of $x + \alpha p$ by $\|x + \alpha p\|$.
In particular, division is the most expensive among the four arithmetic operations.
Therefore, reducing the number of computations of normalization is expected to reduce the total time.
To examine the computational cost reduction, in the following experiments, we deal with the Rayleigh quotient minimization problem on $S^{n-1}$, that is, 
\begin{equation}
\label{prob:sphere}
\begin{array}{ll}
\displaystyle \Mini_{x \in S^{n-1}} \ & f(x) = x^\top Ax,
\end{array}
\end{equation}
where $A \in \mathbb{R}^n$ is a constant symmetric matrix.
The Euclidean gradient of $f \colon \mathbb{R}^n \to \mathbb{R} \colon x \mapsto x^\top A x$ is equal to $2Ax$, and therefore, the Riemannian gradient of $f \colon S^{n-1} \to \mathbb{R} \colon x \mapsto x^\top A x$ is 
\begin{equation}
\label{eq:sphereGrad}
\grad f(x) = P_x(2Ax) = 2(I_n - xx^\top)Ax,
\end{equation}
where $P_x \colon \mathbb{R}^n \to T_x S^{n-1} \colon d \mapsto (I_n - xx^\top)d$ is the orthogonal projection.
Furthermore, the Riemannian Hessian is computed as follows~\cite{AbsMahSep2008}:
\begin{equation}
\label{eq:sphereHess}
\Hess f(x)[p] = 2((I_n-xx^\top)A-(x^\top Ax)I_n)p.
\end{equation}

On the Stiefel manifold $\St(r,n)$, a practical retraction is based on the QR decomposition, which is defined for $X \in \St(r,n)$ and $Y \in T_X \! \St(r,n) = \{Z \in \mathbb{R}^{n \times r} \mid X^\top Z + Z^\top X = I_r\}$ as
\begin{equation}
R_X(Y) = \qf(X+Y),
\end{equation}
where $\qf(W)$ for a full-rank matrix $W \in \mathbb{R}^{n \times r}$ is the Q-factor of the QR decomposition of $W$, that is, if $W = QR$ with $Q \in \mathbb{R}^{n \times r}$ satisfying $Q^\top Q = I_r$ and $R \in \mathbb{R}^{p \times p}$ being an upper triangular matrix with positive diagonal elements (such a decomposition is unique), then $\qf(W) = Q$.
In addition to just computing $X + \alpha Y$, by definition of the retraction, computing $R_X(\alpha Y) = \qf(X + \alpha Y)$ requires extra costs for the QR decomposition of the $n \times r$ matrix $X + \alpha Y$, which is $O(nr^2)$ when $n \gg r$~\cite{golub2012matrix}.
Therefore, at each iterate $k$, when it takes $(l_k+1)$ iterations in the backtracking procedure to find a step length $\alpha_k$, Algorithm~\ref{Armijo_algo} requires $l_k$ times of computation of QR decomposition, each of which costs $O(nr^2)$, whereas Algorithm~\ref{alg:Newton_method} requires only one time of QR decomposition if the Euclidean Armijo condition approximates the Riemannian Armijo condition sufficiently well.
In the following experiments, we deal with the Brockett cost minimization problem on $\St(r,n)$, that is, 
\begin{equation}
\label{prob:Stiefel}
\begin{array}{ll}
\displaystyle \Mini_{X \in \St(r,n)} \ & f(X) = \tr(X^\top AXN),
\end{array}
\end{equation}
where $A \in \mathbb{R}^{n \times n}$ is symmetric, and $N \in \mathbb{R}^{p \times p}$ is a diagonal matrix with diagonal elements $\mu_1 > \mu_2 > \dots > \mu_p > 0$.
Here, $\tr(\cdot)$ denotes the trace of a matrix.
The Euclidean gradient of $f \colon \mathbb{R}^{n \times r} \to \mathbb{R} \colon X \mapsto \tr(X^\top A XN)$ is $2AXN$, and therefore, the Riemannian gradient of $f \colon \St(r,n) \to \mathbb{R} \colon X \mapsto \tr(X^\top A XN)$ is $P_X(2AXN) = 2(AXN - X\sym(X^\top AXN))$, where $P_X \colon \mathbb{R}^{n \times r} \to T_X\! \St(r,n) \colon Y \mapsto Y - X\sym(X^\top Y)$ is the orthogonal projection.
Here, $\sym(W)$ for a square matrix $W \in \mathbb{R}^{n \times n}$ is the symmetric part of $W$ (i.e., $\sym(W) = (W + W^\top) / 2$).

The third example is the SPD manifold $\Sym_{++}(n)$.
On $\Sym_{++}(n)$, we consider the retraction\footnote{This retraction is known as the exponential map for $\Sym_{++}(n)$ equipped with another Riemannian metric $\langle Y, Z\rangle_X \coloneqq \tr(X^{-1}Y X^{-1}Z)$ for $X \in \Sym_{++}(n)$ and $Y, Z \in T_X\!\Sym_{++}(n)$\cite{boumal2023introduction}.} defined for $X \in \Sym_{++}(n)$ and $Y \in T_X \!\Sym_{++}(n) = \Sym(n)$ as
\begin{equation}
R_X(Y) = X^{1/2}\exp(X^{-1/2}YX^{-1/2})X^{1/2} = X\exp(X^{-1}Y).
\end{equation}
Here, although $X\exp(X^{-1}Y)$ is theoretically symmetric, numerically computing $YX^{-1}Y$ may violate the symmetry.
Therefore, for numerically computing the retraction, it should be implemented as
\begin{equation}
R_X(Y) = \sym(X\exp(X^{-1}Y)).
\end{equation}
By definition of the retraction, computing $R_X(Y) = \sym(X\exp(X^{-1}Y))$ requires extra costs, among which computing $X^{-1}Y$ costs $O(n^3)$ and is dominant.
In this case, having computed $X^{-1}Y$, it can be reused to compute $R_X(\alpha Y)$ for $\alpha > 0$ because we have $R_X(\alpha Y) = \sym(X\exp(\alpha(X^{-1}Y)))$.
However, computing the matrix exponential and multiplying the exponential by $X$ from the left remain necessary for each step length candidate.
Therefore, at each iterate $k$, when it takes $l_k$ iterations in backtracking to find an appropriate step length $\alpha_k$, Algorithm~\ref{Armijo_algo} requires the computation of the matrix exponential $(l_k + 1)$ times, that is, $\exp(\beta^l \bar{\alpha}(X_k^{-1}Y_k))$ for $l = 0, 1, \dots, l_k$.
On the other hand, Algorithm~\ref{alg:Newton_method} requires only a single computation of the matrix exponential, that is, $\exp(\beta^{l_k} \bar{\alpha}(X_k^{-1}Y_k))$, if the Euclidean Armijo condition approximates the Riemannian Armijo condition sufficiently well.
Furthermore, since $X$ and $\alpha Y$ are both symmetric, computing the addition $X + \alpha Y$ is both theoretically and numerically symmetric.
Therefore, at this stage, we need not resymmetrize (i.e., use $\sym(\cdot)$ for) $X + \alpha Y$.
In the subsequent experiments, we deal with the following minimization problem on $\Sym_{++}(n)$, namely, 
\begin{equation}
\label{prob:SPD}
\begin{array}{ll}
\displaystyle \Mini_{X \in \Sym_{++}(n)} \ & f(X) = (\det(X)-1)^2.
\end{array}
\end{equation}
Theoretically analyzing the minimizers of this problem is straightforward.
Note that $f(X) \geq 0$ always holds and that, for $I_n \in \Sym_{++}(n)$, $f(I_n) = 0$.
Therefore, the minimum value of $f$ on $\Sym_{++}(n)$ is $0$.
Hence, for $X^{\star} \in \Sym_{++}(n)$, we have
\begin{equation}
\text{$X^{\star}$ is a minimizer of $f$ $\iff$ $f(X^{\star}) = 0$ $\iff$ $\det(X^{\star}) = 1$.}
\end{equation}
The Euclidean gradient of $f \colon \mathbb{R}^{n \times n} \to \mathbb{R} \colon X \mapsto (\det(X)-1)^2$ is equal to $2\det(X)(\det(X)-1)X$, and therefore, the Riemannian gradient of $f \colon \Sym_{++}(n) \to \mathbb{R} \colon X \mapsto (\det(X)-1)^2$ is $\sym(2\det(X)(\det(X)-1)X) = 2\det(X)(\det(X)-1)X$.

\begin{remark}
The matrix space $\mathbb{R}^{m \times n}$ with $m, n \in \mathbb{N}$ can be naturally identified with the Euclidean space $\mathbb{R}^{mn}$ by identifying $X \in \mathbb{R}^{m \times n}$ with $\vect(X) \in \mathbb{R}^{mn}$, where $\vect(X) = \begin{bmatrix}x_{11} & x_{21} & \cdots & x_{m1} & x_{12} & x_{22} & \cdots &x_{mn}\end{bmatrix}$ for $X = [x_{ij}]$.
Then, the standard inner product of $X, Y \in \mathbb{R}^{m \times n}$ is
\begin{equation}
\langle X, Y\rangle = \langle \vect(X), \vect(Y)\rangle = \sum_{i = 1}^m\sum_{j = 1}^n x_{ij}y_{ij} = \tr(X^\top Y).
\end{equation}
Thus, the Stiefel manifold $\St(r,n)$ and the SPD manifold $\Sym_{++}(n)$ are Riemannian submanifolds of the Euclidean spaces $\mathbb{R}^{n \times r}$ and $\mathbb{R}^{n \times n}$, respectively.
\end{remark}

The aforementioned costs for computing addition and retraction are summarized in Table~\ref{tab:costs}.
As listed in the table, computing a retraction requires extra costs compared with computing only the addition.
In other words, the proposed method can reduce such costs by avoiding computing retractions many times.

\begin{table}
\centering
\begin{tabular}{ccc}
Manifold & Costs for $x_k + \alpha_k p_k$ & Extra costs for $R_{x_k}(\alpha_k p_k)$ \\
$S^{n-1} $& $O(n)$ & normalization of $n$-vector\\
$\St(r,n)$ & $O(nr)$ & QR decomp. of $n \times r$ matrix ($O(nr^2)$)\\
$\Sym_{++}(n)$ & $O(n^2)$ & $O(n^3)$ \\
\end{tabular}
\caption{Costs for computing only the addition $x_k + \alpha_k p_k$ and computing a retraction. Note: For $\Sym_{++}(n)$, the retraction considered does not require the addition.}
\label{tab:costs}
\end{table}

\subsection{Numerical results}
In this subsection, we examine the steepest descent method with the proposed line search strategy, namely, Algorithm~\ref{alg:Newton_method} with $H_k = \id_{T_{x_k} \cal M}$.

To ensure the effectiveness of the proposed strategy for different instances and to observe how its performance is affected by the problem size, we applied the conventional Riemannian steepest descent method (with the step length obtained by Algorithm~\ref{Armijo_algo}) and the proposed steepest descent method (Algorithm~\ref{alg:Newton_method} with $H_k = \id_{T_{x_k}\cal M}$) to several instances of the three problems introduced in the previous subsection.

\begin{figure}
\centering
\includegraphics[width=0.5\linewidth]{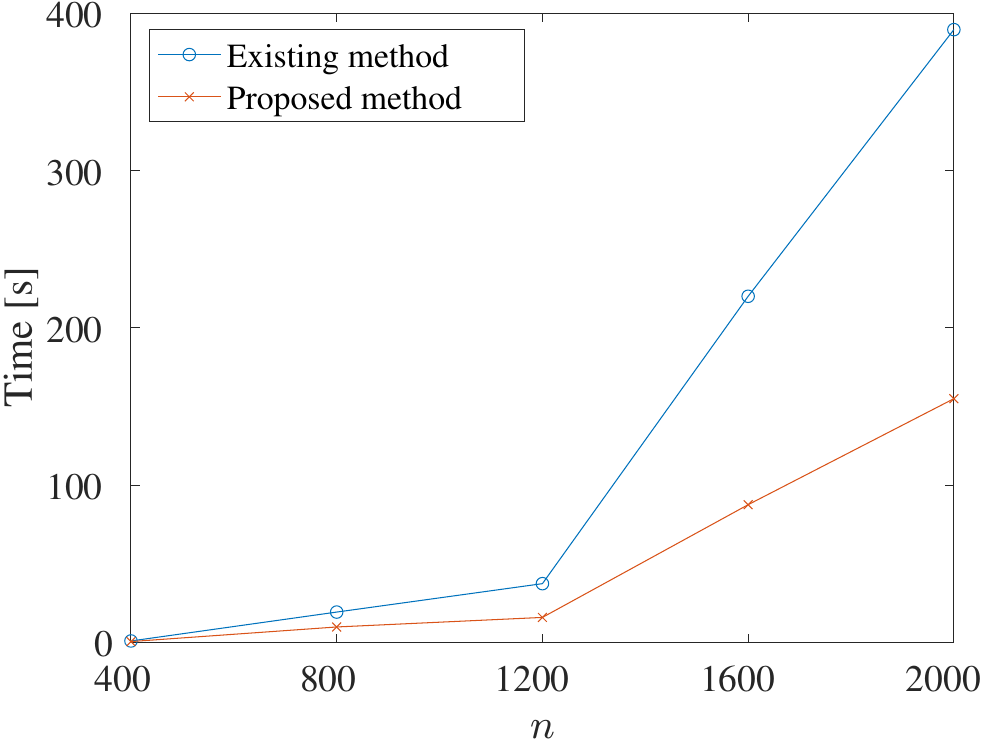}
    \caption{Time comparison of the steepest descent methods with the existing and proposed line search strategies for Problem~\eqref{prob:sphere} on $S^{n-1}$.}
    \label{fig:sphere}
\end{figure}

We first solved the Rayleigh quotient minimization problem~\eqref{prob:sphere} on $S^{n-1}$ with $n = 400, 800, \dots, 2000$.
The initial point $x_0 \in S^{n-1}$ was randomly constructed by \texttt{spherefactory} in Manopt, and the iterations were terminated when $\|\grad f(x_k)\|_{x_k} < 10^{-4}$ was reached.
The result is shown in Figure~\ref{fig:sphere}, which illustrates that the proposed method is always faster than the existing one.
Particularly for larger $n$, the proposed method significantly reduces the computational time.
Table~\ref{tab:sphere} lists the detailed results.
From Table~\ref{tab:sphere}, we can observe that the computational time regarding the proposed method is always less than that for the existing method.
Furthermore, in this case, for each $n$, the number of outer iterations until termination (\#Iter.) for the existing and proposed methods coincided.
Therefore, the time difference between the existing and proposed methods stems from the number of evaluations of the retraction mapping (\#Retr. Eval.).
The results presented in Table~\ref{tab:sphere} indicate that the proposed method significantly decreased \#Retr. Eval.
We may further elaborate on this observation.
For each $n$, let $K$ be the number of iterations, e.g., $K = 2230$ for $n = 400$.
Then, \#BT indicates $\sum_{k = 0}^{K-1}l_k$.
As Algorithm~\ref{Armijo_algo} indicates, for the existing method, \#Retr. Eval. always satisfies
\begin{equation}
\text{\#Retr. Eval.} = K + \text{\#BT} = K + \sum_{k = 0}^{K-1}l_k = \sum_{k = 0}^{K-1}(l_k + 1).
\end{equation}
This is trivial because, for each outer iteration number $k$, the backtracking procedure is iterated $l_k$ times, resulting in the computation of $l_k + 1$ retractions.
By contrast, for the proposed method, if the Euclidean Armijo condition always approximates the Riemannian Armijo condition sufficiently well, then a step length satisfying the Euclidean Armijo condition satisfies the Riemannian one.
As a result, the number of retraction evaluations would always be $1$ for each $k$, and \#Retr. Eval. would be equal to \#Iter.
However, this is an extreme ideal situation and does not apply in this case.
Nevertheless, for the proposed method, \#Retr. Eval. is close to \#Iter, as listed in Table~\ref{tab:sphere}.
This contributed to the acceleration of the proposed algorithm.

\begin{table}
\centering
\begin{tabular}{lrrrrrr}
\toprule 
Method & $n$ & Dim. & Time [s] & \#Iter. & \#BT & \#Retr. Eval. \\ 
\midrule
Existing & \multirow{2}{*}{$400$} & \multirow{2}{*}{$399$} & $0.884$ & $2230$ & $14804$ & $17034$\\
Proposed &  &  & $0.509$ & $2230$ & $14804$ & $2375$\\
\midrule
Existing & \multirow{2}{*}{$800$} & \multirow{2}{*}{$799$} & $19.26$ & $10937$ & $83297$ & $94234$\\
Proposed &  &  & $9.79$ & $10937$ & $83297$ & $11436$\\
\midrule
Existing & \multirow{2}{*}{$1200$} & \multirow{2}{*}{$1199$} & $37.32$ & $9352$ & $73304$ & $82656$\\
Proposed & &  & $15.86$ & $9352$ & $73304$ & $10902$\\
\midrule
Existing & \multirow{2}{*}{$1600$} & \multirow{2}{*}{$1599$} & $220.0$ & $25128$ & $216131$ & $241259$\\
Proposed &  &  & $87.5$ & $25128$ & $216131$ & $25944$\\
\midrule
Existing & \multirow{2}{*}{$2000$} & \multirow{2}{*}{$1999$} & $389.6$ & $33395$ & $300229$ & $333624$\\
Proposed &  &  & $154.9$ & $33395$ & $300229$ & $34473$\\ \bottomrule
\end{tabular}
\caption{Results of the steepest descent method. Problem size $n$, dimension of the manifold (Dim.), computational time (Time), the number of outer iterations (\#Iter.), total of the numbers of backtracking iterations (\#BT), and the number of retraction evaluations (\#Retr. Eval.) until termination for the Rayleigh quotient minimization problem~\eqref{prob:sphere} on $S^{n-1}$ with $n = 400, 800, \dots, 2000$.}
\label{tab:sphere}
\end{table}

We proceed to the second problem~\eqref{prob:Stiefel} on the Stiefel manifold $\St(r,n)$ with $(n, r) = (20, 5), (40, 10), \dots, (100, 25)$.
A similar experiment was conducted; the results are summarized in Table~\ref{tab:Stiefel}.
The initial point $X_0 \in \St(r,n)$ was randomly constructed by \texttt{stiefelfactory} in Manopt, and the iterations were terminated when $\|\grad f(X_k)\|_{X_k} < 10^{-4}$.
For this problem, the numbers of outer iterations of the existing and proposed methods were different.

In certain cases, the proposed method required more iterations to terminate than the existing one did, and vice versa.
This implies that, at an outer iterate number $k$, the backtracking procedure required different times of iterations $l_k$ for the proposed and existing methods, leading to different values of step length $\alpha_k$ and different subsequent iterations.
Nevertheless, in the proposed method, the effect of the reduction of \#Retr. Eval. is high, resulting in considerably less computational time for all instances.

\begin{table}
\centering
\begin{tabular}{lrrrrrr}
\toprule
Method & $(n,r)$ & Dim. & Time [s] & \#Iter. & \#BT & \#Retr. Eval. \\
\midrule
Existing & \multirow{2}{*}{$(20, 5)$} & \multirow{2}{*}{$90$} & $0.471$ & $1853$ & $8857$ & $10710$\\
Proposed &  &  & $0.128$ & $1909$ & $9163$ & $2272$\\
\midrule
Existing & \multirow{2}{*}{$(40, 10)$} & \multirow{2}{*}{$355$} & $36.27$ & $104505$ & $710154$ & $814659$\\
Proposed &  &  & $6.64$ & $95570$ & $649441$ & $110551$\\
\midrule
Existing & \multirow{2}{*}{$(60, 15)$} & \multirow{2}{*}{$795$} & $44.18$ & $89876$ & $715601$ & $805477$\\
Proposed &  &  & $10.79$ & $95306$ & $758997$ & $107785$\\
\midrule
Existing & \multirow{2}{*}{$(80, 20)$} & \multirow{2}{*}{$1410$} & $133.3$ & $203275$ & $1772243$ & $1975518$\\
Proposed &  &  & $43.3$ & $238057$ & $2075691$ & $260405$\\
\midrule
Existing & \multirow{2}{*}{$(100, 25)$} & \multirow{2}{*}{$2200$} & $614.9$ & $654074$ & $6045340$ & $6699414$\\
Proposed &  &  & $202.3$ & $646705$ & $5977288$ & $701794$\\
\bottomrule
\end{tabular}
\caption{Results of the steepest descent method. Problem size $(n, r)$, dimension of the manifold (Dim.), computational time (Time), the number of outer iterations (\#Iter.), a total of the numbers of backtracking iterations (\#BT), and the number of retraction evaluations (\#Retr. Eval.) until termination for the Brockett cost minimization problem~\eqref{prob:Stiefel} on $\St(r,n)$ with $(n, r) = (20, 5), (40, 10), \dots, (100, 25)$.}
\label{tab:Stiefel}
\end{table}

The third problem~\eqref{prob:SPD} was solved on the SPD manifold $\Sym_{++}(n)$ with $n = 200, 400, \dots, 1000$.
To avoid tremendously large values of $\det(X_0)$ and subsequent $\det(X_k)$, we set the initial point $X_0$ as $X_0 \coloneqq I_n + \sym(X_0') / 1000$, which is close to an optimal solution $I_n$, where each element of $X_0'$ was chosen from the interval $[-0.5, 0.5]$ uniformly at random.
The iterations were again terminated when $\|\grad f(X_k)\|_{X_k} < 10^{-4}$.
\begin{table}
\centering
\begin{tabular}{lrrrrrr}
\toprule
Method & $n$ & Dim. & Time [s] & \#Iter. & \#BT & \#Retr. Eval. \\
\midrule
Existing & \multirow{2}{*}{200} & \multirow{2}{*}{20100} & 0.541 & 13 & 104 & 117\\
Proposed &  &  & 0.232 & 13 & 104 & 13\\
\midrule
Existing & \multirow{2}{*}{400} & \multirow{2}{*}{80200} & 2.855 & 14 & 126 & 140\\
Proposed &  &  & 1.387 & 14 & 126 & 14\\
\midrule
Existing & \multirow{2}{*}{600} & \multirow{2}{*}{180300} & 2.261 & 5 & 50 & 55\\
Proposed &  &  & 0.951 & 5 & 50 & 5\\
\midrule
Existing & \multirow{2}{*}{800} & \multirow{2}{*}{320400} & 10.87 & 15 & 150 & 165\\
Proposed &  &  & 4.12 & 15 & 150 & 15\\
\midrule
Existing & \multirow{2}{*}{1000} & \multirow{2}{*}{500500} & 257.2 & 189 & 1890 & 2079\\
Proposed &  &  & 87.3 & 189 & 1890 & 189\\
\bottomrule
\end{tabular}
\caption{Results of the steepest descent method. Problem size $n$, Dimension of the manifold (Dim.), Computational time (Time), the number of outer iterations (\#Iter.), total of the numbers of backtracking iterations (\#BT), and the number of retraction evaluations (\#Retr. Eval.) until termination for the minimization problem \eqref{prob:SPD} on $\Sym_{++}(n)$ with $n = 200, 400, \dots, 1000$.}
\label{tab:SPD}
\end{table}
Table~\ref{tab:SPD} also shows the superiority of the proposed method.
In this case, similar to the first example, \#Iter. for the existing and proposed methods coincided for each value of $n$.
Furthermore, \#Retr. Eval. is always equal to \#Iter for the proposed method.
This means that, for this case, the computed step length candidate satisfying the Euclidean Armijo condition always satisfied the Riemannian Armijo condition.
Therefore, especially for larger $n$, the proposed method significantly reduced \#Retr. Eval., and hence, the computational time.

\begin{remark}
We also performed some preliminary experiments for the Newton method with the proposed line search strategy, that is, Algorithm~\ref{alg:Newton_method} with $H_k = \Hess f(x_k)$.
We applied the Newton method with the existing and proposed line strategies to the Rayleigh quotient minimization problems~\eqref{prob:sphere} on $S^{n-1}$ with several values of $n$.
To guarantee the positive definiteness of the Hessian of $f$, we chose initial points that are close to optimal solutions.
The results showed that, for each $n$, the indicators \#Iter., \#BT, and \#Retr. Eval. always coincided between the existing and proposed methods, with \#BT being always $0$.
This means that the first trials for step length at each iteration always satisfied both the Euclidean and Riemannian Armijo conditions.
Furthermore, the computational time for both methods was very similar.
In fact, in this case, the proposed method required the computation of one more addition $x_k + \alpha p_k$ at each iteration in Step 2 of Algorithm~\ref{alg:Newton_method}.
However, this cost is extremely low and has a slight impact on the overall computational time.
Hence, we conclude that the proposed approach is more effective for the first-order optimization methods than for the second-order ones.
\end{remark}

\section{Concluding remarks} \label{sec:Conclusion}
In this study, a novel line search method that improves the conventional Riemannian Armijo line search was proposed for solving optimization problems on Riemannian submanifolds of the Euclidean spaces.
Existing methods necessitate the computation of the retraction regarding the search direction at each iteration of the backtracking procedure; in contrast, the proposed line search decreases the computational cost by computing a retraction only when necessary.
Moreover, a class of Riemannian optimization methods, including the steepest descent and Newton methods, with the novel line search strategy was proposed, and its global convergence was proved.
Furthermore, we conducted numerical experiments, verifying the practical effectiveness of the proposed method.

\section*{Acknowledgments}
This study was partly supported by grant numbers JP20K14359 and JP25K07125 (Sato), JP21K17709 (Yamakawa), and JP21K11925 and JP24K14985 (Aihara) from the Grants-in-Aid for Scientific Research Program (KAKENHI) of the Japan Society for the Promotion of Science (JSPS).

\section*{Conflict of interest}
The authors declare that they have no conflicts of interest.

\bibliographystyle{abbrv}
\bibliography{sato_bib}
\end{document}